\documentclass[12pt]{article}
\usepackage{t1enc}
\usepackage[latin1]{inputenc}
\usepackage[english]{babel}
\usepackage{amsmath,amsthm}
\usepackage{amsfonts}
\usepackage{amssymb}
\usepackage{enumitem}
\usepackage{latexsym}
\usepackage{graphicx}
\usepackage{epstopdf}
\usepackage{float}
\usepackage{color}
\usepackage[natural]{xcolor}
\usepackage{float,color,fancybox,shapepar,setspace,hyperref}
\theoremstyle{plain}

\newtheorem{theorem}{Theorem}

\newtheorem{claim}{Claim}
\newtheorem{lemma}[claim]{Lemma}
\numberwithin{claim}{section}

\newtheorem*{conjecture*}{Conjecture}
\newtheorem{definition}[claim]{Definition}
\newtheorem{fact}[claim]{Fact}

\numberwithin{remark}{section}

\title{\LARGE On sufficient conditions for planar graphs to be 5-flexible\thanks{This work is supported by NSFC(11971270, 11631014) of China and Shandong Province Natural Science Foundation (ZR2018MA001,ZR2019MA047) of China }}
\author{ Fan Yang\thanks{E-mail address: yangfan5262@163.com.}\\
{\small School of Mathematics, Shandong University, Jinan 250100, China}\\
}
\date{}

\begin{document}
\maketitle

\baselineskip 0.65cm
\begin{abstract}
In this paper, we study the flexibility of two planar graph classes $\mathcal{H}_1$, $\mathcal{H}_2$, where $\mathcal{H}_1$, $\mathcal{H}_2$ denote the set of all hopper-free planar graphs and house-free planar graphs, respectively. Let $G$ be a planar graph with a list assignment $L$. Suppose a preferred color is given for some of the vertices. We prove that if $G\in \mathcal{H}_1$ or $G\in \mathcal{H}_2$ such that all lists have size at least $5$, then there exists an $L$-coloring respecting at least a constant fraction of the preferences.

\baselineskip 0.6cm {\bf Key words:} Planar graph, reducible, discharging, flexibility, satisfiable.
\end{abstract}
\section{Introduction}
\begin{figure}[H]
\begin{center}
\includegraphics[scale=0.8]{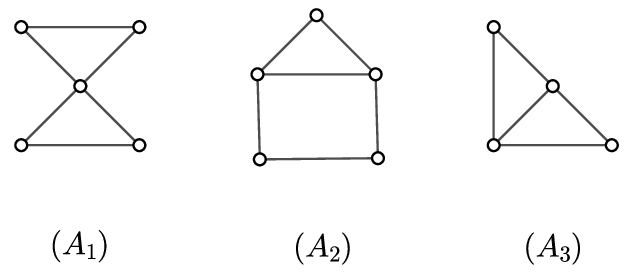}\\
\caption{Special subgraphs}
\label{hhfig}
\end{center}
\end{figure}

\indent
All graphs considered are simple, finite, and loopless, and we follow \cite{Bon} for the terminologies and notation not defined here. Two triangles which intersect exactly at one vertex form a \emph{hopper}, see Figure \ref{hhfig} $(A_1)$. A triangle shares exactly one edge with a $4$-cycle form a \emph{house}, see Figure \ref{hhfig} $(A_2)$. Given a graph $G$, $G$ is called \emph{hopper-free} (or \emph{house-free}) if $G$ does not contain any hopper (or house) as subgraphs. For brevity, denote by $\mathcal{H}_1$ and $\mathcal{H}_2$ the set of all hopper-free planar graphs and house-free planar graphs, respectively. In a proper coloring, we want to assign to each vertex of $G$ one of a fixed number of colors in such a way that adjacent vertices receive distinct colors.
A \emph{list assignment} $L$ for $G$ is a function that assigns to every vertex of $G$ a set (list) $L(v)$ of colors. An \emph{$L$-coloring} is a proper coloring $\phi$ such that $\phi(v)\in L(v)$ for all $v\in V(G)$. If $G$ has a proper coloring $\phi$ such that $\phi(v)\in L(v)$ for each vertex $v$ of $G$, then we say that $G$ is \emph{$L$-colorable}. In addition, we say $L$ is an \emph{$f$-assignment} if $|L(v)|\geq f(v)$ for all $v\in V(H)$. Specifically, $L$ is called a \emph{$k$-assignment} $(k\in \mathbb{N})$ if $f(v)\geq k$ for each $v\in V(G)$. Furthermore, $G$ is \emph{$k$-choosable} if $G$ is $L$-colorable for every $k$-assignment $L$.

Recently, Dvo\v{r}\'{a}k, Norin and Postle introduced a coloring with request as follows. Firstly, we give each vertex of $U\subseteq V(G)$ a preferred color from their list sets, is it possible to properly color $G$ so that at least a constant fraction vertices of $U$ satisfy their preferences?

Initiated by Dvo\v{r}\'{a}k, Norin and Postle \cite{listre}, a \emph{request} for a graph $G$ with a list assignment $L$ is a function $r$ with $\mathrm{dom}(r)\subseteq V(G)$ such that $r(v)\in L(v)$ for all $v\in \mathrm{dom}(r)$. For $\varepsilon >0$, a request $r$ is \emph{$\varepsilon$-satisfiable} if there exists an $L$-coloring $\phi$ of $G$ satisfying $\phi(v)=r(v)$ for at least $\varepsilon |\mathrm{dom}(r)|$ vertices $v\in \mathrm{dom}(r)$. We say that a graph $G$ with the list assignment $L$ is \emph{$\varepsilon$-flexible} if every request is $\varepsilon$-satisfiable. Additionally, we emphasize a stronger weighted form. A \emph{weighted request} is a function $w$ that to each pair $(v,c)$ with $v\in V(G)$ and $c\in L(v)$ assigns a nonnegative real number. Let $w(G,L)=\sum_{v\in V(G), c\in L(v)}w(v,c)$. For $\varepsilon >0$, we say that $w$ is \emph{$\varepsilon$-satisfiable} if there exists an $L$-coloring $\phi$ of $G$ such that
\begin{equation*}
\sum_{v\in V(G)}w(v, \phi(v))\geq \varepsilon w(G,L).
\end{equation*}
We say that $G$ with the list assignment $L$ is \emph{weighted $\varepsilon$-flexible} if every weighted request is $\varepsilon$-satisfiable.

It is worth pointing out that a request $r$ is $1$-satisfiable if and only if the precoloring given by $r$ can be extended to an $L$-coloring of $G$. One can observe that weighted $\varepsilon$-flexibility implies $\varepsilon$-flexibility by giving the request pairs $(v,r(v))$ weight $1$ and all other pairs weight $0$.


Very recently, several scholars contribute a lot on this topic.
For some $\varepsilon>0$, Dvo\v{r}\'{a}k, Norin and Postle \cite{listre} showed that every planar graph is $\varepsilon$-flexible with a $6$-assignment. As we know, planar graphs are $5$-choosable \cite{5cho}, so they conjectured that $6$-assignemnt can be reduced to $5$. In particular, there are lots of results respect to forbidding some configurations in planar graphs. Dvo\v{r}\'{a}k, Masa\v{r}\'{i}k, Mus\'{i}lek and Pangr\'{a}c \cite{DMM} proved that planar graphs without triangles are weighted $\varepsilon$-flexible with a $4$-assignment, the result they gave is the best possible with respect to the list size since planar graphs without triangles are $4$-choosable. Moreover, they also showed that planar graphs of girth at least six are weighted $\varepsilon$-flexible with a $3$-assignment \cite{DMgirth}. However, Thomassen \cite{girth5} studied that planar graphs of girth at least five is $3$-choosable, so there is still a gap left open. Masa\v{r}\'{i}k \cite{MT} showed that $C_4$-free planar graphs are weighted $\varepsilon$-flexible with a $5$-assignment. Since planar graphs without $4$-cycles are $4$-choosable \cite{Xu}, Masa\v{r}\'{i}k conjectured that his result about list assignment would be reduced to $4$. In addition, Choi et.al \cite{CHO} proved three theorems: every planar graph (i) without $K_4^{-}$ is weighted $\varepsilon$-flexible with a $5$-assignment, (ii) without $C_4$ and $C_3$ distance at least $2$ is weighted $\varepsilon$-flexible with a $4$-assignment, (iii) without $C_4$, $C_5$, $C_6$ is weighted $\varepsilon$-flexible with a $4$-assignment. Their first theorem has strengthened the result of Masa\v{r}\'{i}k, which is a good bound up to the list size compared to choosability, since the conjecture that
$K_4^{-}$-free planar graphs are $4$-choosable is still open. Nowadays, Yang and the author \cite{YY} extended the third theorem of Choi et. al, they showed that every $\{C_4,C_5\}$-free planar graph is weighted $\varepsilon$-flexible with a $4$-assignment, which is the best possible with respect to the list size, since Voigt \cite{Voi} gave a planar graph without $C_4$ and $C_5$ is not $3$-choosable.

In the paper, we mainly investigate the weighted $\varepsilon$-flexibility of two classes of planar graphs.

\begin{theorem}\label{thm1}
If $G$ is hopper-free, then there exists $\varepsilon >0$ such that $G$ is weighted $\varepsilon$-flexible with a $5$-assignment.
\end{theorem}

\begin{theorem}\label{thm2}
If $G$ is house-free, then there exists $\varepsilon >0$ such that $G$ is weighted $\varepsilon$-flexible with a $5$-assignment.
\end{theorem}

Until now, no result states that $G$ is $4$-choosable if $G\in \mathcal{H}_1$ or $G\in \mathcal{H}_2$. However, Wang and Lih \cite{hop4cho} proved that a planar graph $H$ is $4$-choosable if $H$ has no intersecting $3$-cycle, that is, $H$ contains neither hopper nor diamond (which is the graph isomorphic to $(A_3)$, see Figure \ref{hhfig}). On the other hand, Borodin and Ivanova \cite{hou4cho}, Cheng et al. \cite{CCW} showed that a planar graph $H$ is $4$-choosable if $H$ contains no $4$-cycle which shares at least one common edge with a $3$-cycle, that is, $H$ contains neither house nor diamond.

The rest of the paper is organized as follows. In Section 2, we develop the notation and introduce some essential tools used in list coloring settings. In Section \ref{hpgra}, we give the proof of Theorem \ref{thm1}, the main idea is to produce some reducible configurations and then use discharging method to get a contradiction. In Section \ref{hfgra}, we prove Theorem \ref{thm2}.

\section{Preliminaries}\label{prelim}
We shall split Section \ref{prelim} into two parts. In Section \ref{defini}, we first give some definitions related to our topic. In Section \ref{tool}, we shall present several essential tools to the proof of our theorems.

\subsection{Definitions}\label{defini}
Let \emph{$1_{S}$} denote the characteristic function of $S$, i.e., $1_{S}(v)=1$ if $v\in S$ and $1_{S}(v)=0$ otherwise. For functions that assign integers to vertices of $H$, we define addition and subtraction in the natural way, adding$/$subtracting their values at each vertex independently. Given a graph $H$ and a vertex $v \in V(H)$. For a function $f: V(H)\rightarrow \mathbb{Z}$ and a vertex $v\in V(H)$, let \emph{$f\downarrow v$} denote the function such that $(f\downarrow v)(w)=f(w)$ for $w\neq v$ and $(f\downarrow v)(v)=1$. Given a set of graphs $\mathcal{F}$ and a graph $H$, a set $S\subseteq V(H)$ is $\mathcal{F}$-forbidding if the graph $H$ together with one additional vertex adjacent to all of the vertices in $S$ does not contain any graph from $\mathcal{F}$. We first give a crucial definition of \emph{$(\mathcal{F},k)$-boundary-reducible} as follows.

\begin{definition}
A graph $H$ is an $(\mathcal{F},k)$-boundary-reducible induced subgraph of a graph $G$ if there exists a set $B\subsetneqq V(H)$ such that
\begin{description}
\item[(\textbf{FIX})] for every $v\in V(H)\backslash B$, $H-B$ is $L$-colorable for every $((k-\deg_{G}+\deg_{H-B})\downarrow v)$-assignment $L$, and
\item[(\textbf{FORB})] for every $\mathcal{F}$-forbidding set $S\subseteq V(H)\backslash B$ of size at most $k-2$, $H-B$ is $L'$-colorable for every $(k-\deg_{G}+\deg_{H-B}-1_{S})$-assignment $L'$.
\end{description}
\end{definition}

We will call the set $B$ in (\textbf{FORB}) as the \emph{boundary} of the configuration in the following discussion. By the definition of (\textbf{FORB}), we get that (\textbf{FORB}) is implied by (\textbf{FIX}) when $|S|=1$. Hence in the following discussion, we mainly consider the case $2\leq|S|\leq k-2$.

\begin{definition}
Let $G$ be a graph with lists of size $k$ that does not contain any graph in $\mathcal{F}$ as an induced subgraph. We define $(\mathcal{F},k,b)$-resolution of $G$ as a set $G_i$ of nested subgraphs for $0\leq i\leq M$, such that $G_0:=G$ and
$$G_i:=G-\bigcup_{j=1}^{i}(H_j-B_j),$$
where each $H_i$ is an induced $(\mathcal{F},k)$-boundary-reducible subgraph of $G_{i-1}$ with boundary $B_i$ such that $|V(H_i)\backslash B_i|\leq b$ and $G_M$ is an $(\mathcal{F},k)$-boundary-reducible subgraph with empty boundary and size at most $b$. For technical reasons, let $G_{M+1}:=\emptyset$.
\end{definition}

Our strategy is to prove every graph that does not contain any graph from $\mathcal{F}$ as a subgraph contains a reducible subgraph. Actually, we regard a resolution as an inductively-defined object obtained by iteratively identifying some reducible subgraph $H$ with boundary $B$ and deleting $H-B$ until $V(G)$ is exhausted.

\subsection{Basic tools}\label{tool}
The following lemma derived from Choi et al. provide us with a unified approach to deal with the weighted flexibility of any graph with forbidden subgraphs, which also strengthen the key lemma implicitly presented by Dvo\v{r}\'{a}k, Norin and Postle in \cite{listre}, and explicitly formulated as Lemma 4 in \cite{DMM}.

\begin{lemma}[\cite{CHO}]
For all integers $k\geq3$ and $b\geq1$ and for all sets $\mathcal{F}$ of forbidden subgraphs, there exists an $\varepsilon>0$ as follows. Let $G$ be a graph with an $(\mathcal{F},k,b)$-resolution. Then $G$ with any assignment of lists of size $k$ is weighted $\varepsilon$-flexible.
\end{lemma}

The well-known lemma below provide us with a method to deal with the coloring problem of a graph, which will be used frequently in our proofs.
\begin{lemma}[\cite{Thomasdeg}]\label{degree}
Let $G$ be a connected graph and $L$ a list assignment such that $|L(u)|\geq \deg(u)$ for all $u\in V(G)$. If either there exists a vertex $u\in V(G)$ such that $|L(u)|> \deg(u)$, or some $2$-connected component of $G$ is neither complete nor an odd cycle, then $G$ is $L$-colorable.
\end{lemma}

\section{Proof of Theorem \ref{thm1}}\label{hpgra}
In this section, we shall first collect essential notation and then find some reducible subgraphs. Finally we use Euler's formula to complete the proof of Theorem \ref{thm1}.

\subsection{Notation}\label{nota}
A plane graph is a particular drawing of a planar graph in the Euclidean plane.
Let $G$ be a plane graph, let us denote by $V(G)$, $E(G)$, $F(G)$ the vertex set, edge set, face set of $G$, respectively. We denote by $d(v)$ and $\delta(G)$ the degree of a vertex $v\in V(G)$ and minimum degree of $G$. A vertex $v$ is called a $k$-vertex, a $k^{+}$-vertex or a $k^{-}$-vertex if $d(v)=k$, $d(v)\geq k$ or $d(v)\leq k$, respectively. For any face $f\in F(G)$, the degree of $f$, denoted by $d(f)$, is the length of the shortest boundary walk of $f$, where each cut edge is counted twice. A $k$-face, a $k^+$-face or a $k^{-}$-face is a face of degree $k$, degree at least $k$, or degree at most $k$, respectively. We write $f=(d_1, \ldots  d_n)$ if $v_1, \ldots, v_n$ are the boundary vertices of $f$ with $d(v_i)=d_i$ for all $i\in\{1,2, \ldots, n\}$. We say that $f=(d_1^{+}, \ldots  d_n)$ if $d(v_1)\geq d_1$ and $d(v_i)=d_i$ for all $i\in\{2, \ldots, n\}$; and similarly for other combinations. In addition, Let $f_k(v)$, $n_k(f)$ denote the number of $k$-faces incident with the vertex $v$ and the number of $k$-vertices incident with the face $f$, respectively.




\subsection{Reducible configurations}
\begin{figure}[H]
\begin{center}
\includegraphics[scale=0.8]{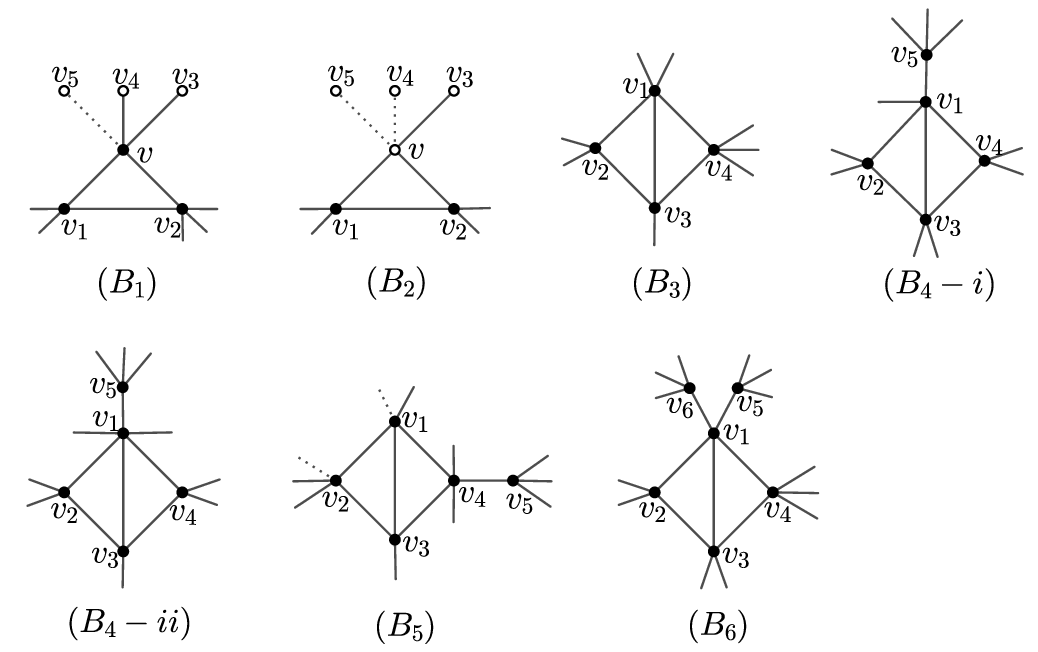}\\
\caption{Reducible subgraphs. The edge marked with a dashed line may not exist, and its existence depends on the degree of the vertex incident with it.}
\label{hopfig}
\end{center}
\end{figure}

Note that in all figures of the paper, any vertex marked with $\bullet$ has no edges of $G$ incident with it other than those shown. In the following, we say a vertex $u$ has $\gamma$ ($\gamma\in \mathbb{N}$) \emph{available} colors in a configuration $H$ means that the maximum number of colors remaining in $L(u)$ is $\gamma$ after coloring vertices exterior to $H$. When considering (\textbf{FIX}), we reduce the number of available colors on the vertex to $1$, and when considering (\textbf{FORB}), we reduce the number of available colors on the vertices in $S$ by $1$.
\begin{lemma}\label{configu1}
Let $G\in \mathcal{H}_1$. If $G$ contains one of the following configurations (see Figure \ref{hopfig}).
\begin{description}
\item[$(B_1)$] A cycle $vv_1v_2$ such that $4\leq d(v)\leq5$, $d(v_i)=4$ for each $i \in \{1,\ldots, d(v)\}\backslash \{2\}$ and $d(v_2)=5$, $v_j \in N(v)$ for each $j\in \{1,\ldots, d(v)\}$.
\item[$(B_2)$] A cycle $vv_1v_2$ such that $4\leq d(v)\leq6$ and $d(v_i)=4$ with $v_i \in N(v)$ for each $i\in \{1,\ldots, d(v)-1\}$.
\item[$(B_3)$] A cycle $v_1v_2v_3v_4$ and an edge $v_1v_3$ such that $d(v_2)=d(v_3)=4$, $d(v_1)=d(v_4)=5$.
\item[$(B_4)$] A cycle $v_1v_2v_3v_4$, an edge $v_1v_3$ and an edge $v_1v_5$ such that one of the following holds,
    \begin{enumerate}
    \item[(i)] $d(v_1)=d(v_3)=5$ and $d(v_2)=d(v_4)=d(v_5)=4$;
    \item[(ii)] $d(v_1)=6$, $d(v_i)=4$ for each $i\in\{2,3,4,5\}$.
    \end{enumerate}
\item[$(B_5)$] A cycle $v_1v_2v_3v_4$, and two edges $v_1v_3$, $v_4v_5$ such that $d(v_i)\leq5$ for each $i\in \{1,2\}$, $d(v_4)=5$, and $d(v_3)=d(v_5)=4$.
\item[$(B_6)$] A cycle $v_1v_2v_3v_4$, and three edges $v_1v_3$, $v_1v_5$, $v_1v_6$ such that $d(v_i)=5$ for each $i\in\{1,3,4\}$, $d(v_j)=4$ for each $j\in\{2,5,6\}$.
\end{description}
Then $G$ contains a $(\mathcal{H}_1,5)$-boundary-reducible induced subgraph with empty boundary.
\end{lemma}

\begin{proof}
Let $H_i$ be the graph isomorphic to one of $(B_i)$ for $i\in \{1,\ldots,6\}$ and set the boundary to be empty. It suffices to prove that $H_i$ satisfies (\textbf{FIX}) and (\textbf{FORB}) for each $i\in \{1,\ldots,6\}$.

For $H_1$. We only consider $d(v)=4$ since the same arguments yield to $d(v)=5$.

(\textbf{FIX}): Note that $v$ has five available colors, $v_1$ has three available colors, and $v_i$ has two available colors for each $i\in\{2,3,4\}$. If $v$ (or $v_1$) is fixed, then we first color $v$ (or $v_1$). Finally $H_1-v$ (or $H_1-v_1$) can be colored by Lemma \ref{degree}. If $v_2$ is fixed, then we can greedily color $v_2$, $v_1$, $v$, $v_3$, $v_4$ in order. Fixing any other vertex in $\{v_3,v_4\}$ is handled in a similar fashion.

(\textbf{FORB}): Let $S\subseteq V(H_1)$ of size at most $3$. Recall that (\textbf{FORB}) is implied by (\textbf{FIX}) when $|S|=1$. So we mainly discuss the case $2\leq |S|\leq3$ in the following proof. Suppose $|S|=2$. Then $S$ can be chosen as the following pairs: $\{v,v_1\}$, $\{v,v_2\}$, $\{v_1,v_2\}$, $\{v_1,v_3\}$, $\{v_1,v_4\}$, $\{v_2,v_3\}$, $\{v_2,v_4\}$, $\{v_3,v_4\}$. If $S=\{v,v_1\}$ (or $S=\{v,v_2\}$), then we can greedily color $v_2$, $v_1$, $v$, $v_3$, $v_4$ in order. If $S=\{v_1,v_3\}$ (or $S=\{v_1,v_4\}$), then $H_1$ can be colored by Lemma \ref{degree}. If $S=\{v_1,v_2\}$, then we can greedily color $v_2$, $v_1$, $v$, $v_3$, $v_4$ in order. Forb the remaining pair of vertices described as above can be handled in a similar fashion, so we omit them. Now we discuss the case $|S|=3$. By the definition of $S$, we know that $S$ can be chosen as the following triples: $\{v, v_1, v_2\}$, $\{v_1, v_2, v_3\}$, $\{v_1, v_2, v_4\}$, $\{v_2, v_3, v_4\}$. If $S=\{v, v_1, v_2\}$, then we greedily color $v_2$, $v_1$, $v$, $v_3$, $v_4$ in order. The rest triples can be handled in a similar fashion.

For $H_2$. We only consider the case $d(v)=5$ as the analysis of the other two cases are similar.

(\textbf{FIX}): Note that $v$ has four available colors, $v_i$ has three available colors for each $i\in\{1,2\}$, and $v_j$ has two available colors for each $j\in\{3,4\}$. If $v$ is fixed, then we greedily color $v$, $v_1$, $v_2$, $v_3$, $v_4$ in order. If $v_i$ is fixed for each $i\in\{1,2,3,4\}$, then we first color $v$, and finally the coloring can be extended to $H_2-v_i$ by Lemma \ref{degree}.

(\textbf{FORB}): When $|S|=2$, if $S=\{v, v_j\}$ $(j\in\{1,2\})$, then we can greedily color $v_j$, $v_{3-j}$, $v$, $v_3$, $v_4$ in order. If $S=\{s_1, s_2\}$ $(s_1\neq s_2)$, where $s_i\in\{v_1,v_2,v_3,v_4\}$ for each $i$. Similarly, we first greedily color the vertices in $S$ and then extend the coloring to the remaining vertices, which is possible since $v$ has four available colors. When $|S|=3$, if $S=\{v, v_1, v_2\}$, then we can color $v_1$, $v_2$, $v$, $v_3$, $v_4$ in order. If $S=\{s_1,s_2,s_3\}$, where $s_i\in\{v_1, v_2, v_3, v_4\}$ for each $i$ and $s_i$ are distinct from each other, then we can first greedily color the vertices in $S$ and finally extend the coloring to the remaining vertices. As a result, (\textbf{FORB}) holds.

For $H_3$. We verify both (\textbf{FIX}) and (\textbf{FORB}) holds.

(\textbf{FIX}): Note that $v_i$ has three available colors for each $i\in \{1,2\}$, $v_3$ has four available colors, and $v_4$ has two available colors. If $v_i$ is fixed for each $i\in\{1,2,3,4\}$, then we first color $v_i$. At last, $H_3-v_i$ can be colored by Lemma \ref{degree}.

(\textbf{FORB}): Note that $|S|\leq2$. If $S=\{v_1, v_3\}$, then we greedily color $v_1$, $v_4$, $v_3$, $v_2$ in order. If $S=\{v_2, v_4\}$, then we greedily color $v_4$, $v_1$, $v_2$, $v_3$ in order.

For $H_4$. Suppose $d(v_1)=d(v_3)=5$ and  $d(v_2)=d(v_4)=d(v_5)=4$.

(\textbf{FIX}): Note that $v_1$ has four available colors, $v_i$ has three available colors for each $i\in \{2,3,4\}$, and $v_5$ has two available colors. If $v_1$ is fixed, then we greedily color $v_1$, $v_5$, $v_2$, $v_3$, $v_4$ in order. If $v_i$ is fixed for each $i\in\{2,3,4,5\}$, then we first color $v_i$. Finally, $H-v_i$ can be colored by Lemma \ref{degree}.

(\textbf{FORB}): We first consider $|S|=2$. If $S=\{v_1, v_3\}$, then we greedily color $v_3$, $v_2$, $v_1$, $v_4$, $v_5$ in order. If $S=\{v_2, v_4\}$, then we greedily color $v_2$, $v_3$, $v_4$, $v_1$, $v_5$ in order. If $S=\{v_j, v_5\}$ for each $j\in\{2,3,4\}$, we may assume $j=2$, then we greedily color $v_5$, $v_2$, $v_1$, $v_3$, $v_4$ in order. One can observe that the case $j=3$ or $j=4$ admits (\textbf{FORB}) by the same arguments. If $|S|=3$, then $S=\{v_2, v_4, v_5\}$, then we first color $v_5$, and then $H-v_5$ can be colored by Lemma \ref{degree}, (\textbf{FORB}) holds.

Suppose $d(v_1)=6$, $d(v_i)=4$ for each $i\in\{2,3,4,5\}$.

(\textbf{FIX}): Note that $v_i$ has three available colors for each $i\in \{1,2,4\}$, $v_3$ has four available colors, and $v_5$ has two available colors. Fix any other vertex, say $v_1$, then we can first color $v_1$ and then extend the coloring to the remaining vertices.

(\textbf{FORB}): When $|S|=2$, if $S=\{v_1, v_3\}$, then we color $v_5$, $v_1$, $v_3$, $v_2$, $v_4$ in order. If $S=\{v_2, v_4\}$, then we greedily color $v_2$, $v_4$, $v_1$, $v_5$, $v_3$ in order. If $S=\{v_j, v_5\}$ for $j\in\{2,3,4\}$, say $j=2$, then we greedily color $v_5$, $v_2$, $v_1$, $v_4$, $v_3$ in order. It is easy to check the cases $j=3$ and $j=4$ also admit (\textbf{FORB}) by the same arguments. If $|S|=3$, then $S=\{v_2, v_4, v_5\}$, and thus we can greedily color $v_5$, $v_4$, $v_1$, $v_2$, $v_3$ in order, implying (\textbf{FORB}).

For $H_5$. We may assume that $d(v_1)=d(v_2)=5$.

(\textbf{FIX}): Note that $v_i$ has three available colors for each $i\in \{1,4\}$, $v_j$ has two available colors for each $j\in \{2,5\}$, and $v_3$ has four available colors. If $v_i$ for some $i\in\{1,2\}$ is fixed, then we greedily color $v_i$, $v_{3-i}$, $v_3$, $v_4$, $v_5$ in order. If $v_4$ is fixed, then we greedily color $v_4$, $v_5$, $v_1$, $v_2$, $v_3$ in order. If $v_i$ is fixed for some $i\in\{3,5\}$, then we first color $v_i$, and then by Lemma \ref{degree}, $H_5-v_i$ can be colored.

(\textbf{FORB}): When $|S|=2$. If $S=\{v_1, v_3\}$, then we greedily color $v_1$, $v_2$, $v_3$, $v_4$, $v_5$ in order. If $S=\{v_2, v_4\}$, then we greedily color $v_2$, $v_1$, $v_4$, $v_5$, $v_3$ in order. If $S=\{v_j, v_5\}$ $(j\in\{1,2,3\})$, we consider $j=1$ here and then we greedily color $v_5$, $v_4$, $v_1$, $v_2$, $v_3$ in order. For $j=2$ or $j=3$, it is easy to check $H_5$ can be colored by the same arguments. When $|S|=3$, we know that $S=\{v_1, v_3, v_5\}$. Let $L'$ be a $(5-\deg_{G}+\deg_{H}-1_{S})$-assignment. We get that $|L'(v_1)|=|L'(v_2)|=2$, $|L'(v_3)|=|L'(v_4)|=3$, and $|L'(v_5)|=1$. If $L'(v_4)\cap L'(v_5)=\emptyset$, then we can $L'$-color $v_5$, $v_1$, $v_2$, $v_3$, $v_4$ in order. Otherwise, we first color $v_5$. Next let $H=H_5\backslash \{v_5\}$ and $L^{*}$ be a assignment for $H$ obtained by $L'$ be removing the only color in $L'(v_5)$ from the list of the the vertex set $\{y:|yv_5\in E(H_5)\}$, that is, $|L^{*}(v_i)|=|L'(v_i)|$ for each $i\in \{1,2,3\}$, $|L^{*}(v_4)|=2$. Now we discuss whether $L^{*}(v_2)\cap L^{*}(v_4)=\emptyset$. If $L^{*}(v_2)\cap L^{*}(v_4)\neq \emptyset$, let $c_1\in L^{*}(v_2)\cap L^{*}(v_4)$, then we first color $v_2$ and $v_4$ with $c_1$ and then greedily $L^{*}$-color $v_1$, $v_3$ in order. Otherwise $L^{*}(v_2)\cap L^{*}(v_4)=\emptyset$, then there must be a color in $L^{*}(v_2)\cup L^{*}(v_4)$ but not in $L^{*}(v_3)$, we denote the color by $c_2$ and assume that $c_2\in L^{*}(v_4)$, then we color $v_4$ with $c_2$ and greedily $L^{*}$-color $v_1$, $v_2$, $v_3$ in order, implying (\textbf{FORB}).

For $H_6$. We shall prove that $H_6$ also satisfies (\textbf{FIX}) and (\textbf{FORB}).

(\textbf{FIX}): Note that $v_1$ has five available colors, $v_i$ has two available colors for each $i\in\{4,5,6\}$, and $v_j$ has three available colors for each $j\in\{2,3\}$. If $v_1$ is fixed, then we greedily color $v_1$, $v_5$, $v_6$, $v_4$, $v_3$, $v_2$ in order. If $v_3$ is fixed, then we greedily color $v_3$, $v_4$, $v_2$, $v_1$, $v_5$, $v_6$ in order. If $v_4$ is fixed, then we greedily color $v_4$, $v_3$, $v_2$, $v_1$, $v_5$, $v_6$ in order. If $v_i$ is fixed for each $i\in\{2,5,6\}$, then we first color $v_i$, and $H-v_i$ can be colored by Lemma \ref{degree}.

(\textbf{FORB}): When $|S|=2$. If $S=\{v_1, v_3\}$, then we can color $v_3$, $v_4$, $v_1$, $v_5$, $v_6$, $v_2$ in order. If $S=\{v_2, v_4\}$, then we can color $v_4$, $v_3$, $v_2$, $v_1$, $v_5$, $v_6$ in order. If $S=\{v_5,v_6\}$, then we can color $v_5$, $v_6$, $v_1$, $v_4$, $v_3$, $v_2$ in order. If $S=\{s_1, s_2\}$, where $s_1\in\{v_2,v_3,v_4\}$, $s_2\in\{v_5,v_6\}$, then we first greedily color $s_1$, $s_2$ and finally extend the coloring to the remaining vertices. When $|S|=3$, then we have $S=\{v_2,v_4,v_i\}$ $(i\in\{5,6\})$, then we can color $v_i$, $v_4$, $v_3$, $v_2$, $v_1$, $v_{11-i}$ in order, implying (\textbf{FORB}).

From all the above cases, both (\textbf{FIX}) and (\textbf{FORB}) hold, and thus $H_i$ is $(\mathcal{H}_1,5)$-boundary-reducible for each $i\in\{1,\ldots,6\}$.
\end{proof}
\subsection{Discharging}
Let $\mathcal{F}=\mathcal{H}_1$, and let $G_1$ be a counterexample to Theorem \ref{thm1} with minimum number of vertices. Fix a plane embedding of $G_1$, by minimality of $G_1$, we get that $G_1$ is connected. Let $L$ be a list assignment on $V(G_1)$ where each vertex receives at least five colors. Note that $G_1$ does not contain any configurations shown in Lemma \ref{configu1}. By \cite{MT}, we get that $G_1$ has no $3^-$-vertex. Since $G_1$ is also a plane graph, by Euler's Formula, we obtain

$$\sum_{v\in V(G_1)}(d(v)-4)+\sum_{f\in F(G_1)}(d(f)-4)=-8. \eqno{(\mathrm{I})}$$

We define an \emph{initial} charge $c$ on $V(G_1)\cup F(G_1)$ by letting
\begin{align*}
\begin{split}
c(x)=\left\{
 \begin{array}{ll}
  d_{G_1}(x)-4    & $if$ \ x=v\in V(G_1), \\
  d_{G_1}(x)-4    & $if$ \ x=f\in F(G_1).
 \end{array}
 \right.
 \end{split}
\end{align*}
We will obtain a \emph{final} charge $\tilde{c}$ from $c$ by discharging rules R1-R6 below. Since these rules merely move charges around, $\mathrm{(I)}$ gives
 $$\sum\limits_{x\in V(G_1)\cup F(G_1)}\tilde{c}(x)=\sum\limits_{x\in V(G_1)\cup F(G_1)}c(x)< 0. \eqno{\mathrm{(II)}}$$
We will get a contradiction by proving $\tilde{c}(x)\geq0$ for each element $x\in V(G_1)\cup F(G_1)$. Since $G_1\in \mathcal{H}_1$, we immediately have the following fact.
\begin{fact}
For each $v\in V(G_1)$ with $d(v)\geq4$, we have $f_3(v)\leq 2$.
\end{fact}

For brevity, $(4,4,5^{-})$-face is called a \emph{bad} face and a $4$-vertex lying on a bad face is called a \emph{bad} $4$-vertex. In addition, a $3$-face $f$ is called a \emph{singleton} if all faces incident with it are $4^{+}$-faces, while two consecutive $3$-faces form a \emph{doubleton} $\bar{f}$. The discharging rules are as follows.

\begin{description}
\item[$\bf R1.$] Each $5^{+}$-vertex $v$ sends $\frac{1}{6}$ to its adjacent bad $4$-vertex which is not lying on the same $3$-face with $v$.
\item[$\bf R2.$] Each bad $4$-vertex sends the total charge it received to its incident bad $3$-face.
\item[$\bf R3.$] Each $5$-vertex $v$ with $f_3(v)=1$ sends $a$ to its incident $3$-face $f$. let $n_b^{*}(v)$ be the number of bad $4$-vertices incident with $v$ while not lying on $f$.
\begin{enumerate}
\item[$\bf R3.1.$] If $f=(4,5^{-},5)$, then
$a=\begin{cases}
1 \ \ & $if$ \ n_b^{*}(v)=0,\\
\frac{2}{3}, \ \ & $if$ \ 1\leq n_b^{*}(v)\leq2.
\end{cases}$
\item[$\bf R3.2.$] If $f=(4,6^{+},5)$, then $a=\frac{1}{2}$.
\item[$\bf R3.3.$] If $f=(5^{+},5^{+},5)$, then
 $a=\begin{cases}
1 \ \ & $if$ \ n_b^{*}(v)=0,\\
\frac{1}{2} \ \ & $if$ \ 1\leq n_b^{*}(v)\leq3.
\end{cases}$
\end{enumerate}
\item[$\bf R4.$] Each $5$-vertex $v$ with $f_3(v)=2$ sends $a$ to its incident doubleton $\bar{f}$, let $n_b(v)$ be the number of bad $4$-vertices not lying on $\bar{f}$, then $a=1-\frac{n_b(v)}{6}$.

\item[$\bf R5.$] Each $6$-vertex $v$ sends $a$ to its incident $3$-face,
\begin{enumerate}
\item[$\bf R5.1.$] If $f_3(v)=1$, then $a=\frac{4}{3}$.
\item[$\bf R5.2.$] If $v$ is incident with a doubleton $\bar{f}$ and there are at most two $4$-vertices on the $\bar{f}$, then $a=\frac{3}{2}$.
\item[$\bf R5.3.$] If $v$ is incident with a doubleton $\bar{f}$ and there are three $4$-vertices on the $\bar{f}$, then $a=2$.
\end{enumerate}
\item[$\bf R6.$] Each $7^{+}$-vertex $v$ sends $a$ to its incident $3$-face,
\begin{enumerate}
\item[$\bf R6.1.$] If $f_3(v)=1$, then $a=\frac{4}{3}$.
\item[$\bf R6.2.$] If $v$ is incident with a doubleton $\bar{f}$, then $v$ sends $2$ to $\bar{f}$.
\end{enumerate}
\end{description}

Let $f$ be a face of $G_1$. If $d(f)\geq4$, the initial charge is not changed, and thus $\tilde{c}(f)=c(f)\geq0$. It remains to consider the case $d(f)=3$. In particular, if two $3$-faces $f_1$ and $f_2$ are consecutive, that is, $f_1$ and $f_2$ form a doubleton $\bar{f}$, In this situation, $c(\bar{f})=c(f_1)+c(f_2)=-2$, and then we discuss the final charge $\tilde{c}(\bar{f})$.

\textbf{Case 1:} $f$ is a singleton .

If $f$ is bad, then by $(B_1)$ and $(B_2)$, each $4$-vertex on the $f$ must be adjacent to two $5^{+}$-vertices, and thus $\tilde{c}(f)\geq -1+\min\{3\times2\times \frac{1}{6}, 2\times \frac{1}{6}+\frac{2}{3}\}=0$ by R1-R3. If $f=(4,4,6^{+})$, then $\tilde{c}(f)\geq -1+\frac{4}{3}=\frac{1}{3}>0$ by R5.1. If $f=(4, 5^{+}, 5^{+})$, then $\tilde{c}(f)\geq -1+\min\{2\times \frac{2}{3}, \frac{1}{2}+\frac{4}{3}\}=\frac{1}{3}>0$ by R3. If $f=(5^{+}, 5^{+}, 5^{+})$, then $\tilde{c}(f)\geq -1+3\times \frac{1}{2}=\frac{1}{2}>0$ by R3.3.

\begin{figure}[H]
\begin{center}
\includegraphics[scale=0.75]{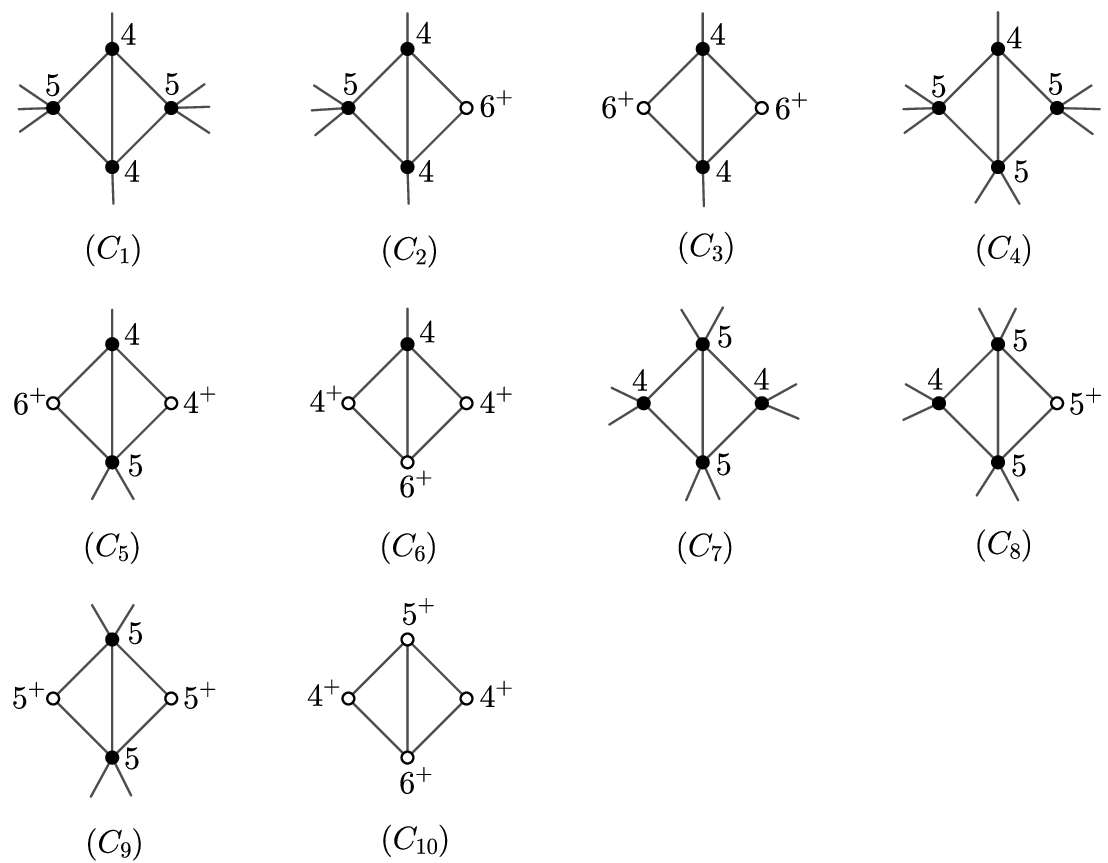}\\
\caption{Different kinds of doubletons, the number near a vertex denotes its degree.}
\label{doub}
\end{center}
\end{figure}

\textbf{Case 2:} $\bar{f}$ is a doubleton  (see Figure \ref{doub}).

It follows from Lemma \ref{configu1} $(B_5)$ that two $5$-vertices in $(C_1)$ are not adjacent to a $4$-vertex any more, thus $\tilde{c}(\bar{f})\geq -2+2\times1=0$ by R3.1. As for $(C_2)$, $\tilde{c}(\bar{f})\geq -2+\frac{2}{3}+\frac{4}{3}=0$ by R3.1, R5.1 and R6.1. For $(C_3)$, we have $\tilde{c}(\bar{f})\geq -2+2\times \frac{4}{3}=\frac{2}{3}>0$ by R5.1 and R6.1. For $(C_4)$, we have $\tilde{c}(\bar{f})\geq -2+3\times\frac{2}{3}=0$ by R3 and R4. For $(C_5)$, $\tilde{c}(\bar{f})\geq -2+\frac{4}{3}+\frac{2}{3}=0$ by R4-R6. For $(C_6)$, if there are three $4$-vertices lying on $\bar{f}$, then $\tilde{c}(\bar{f})\geq -2+2=0$ by R5.3 and R6.2. Otherwise $\tilde{c}(\bar{f})\geq -2+\min\{\frac{1}{2}+\frac{3}{2}, 2\}=0$ by R3, R5 and R6. For $(C_7)$, it follows from Lemma \ref{configu1} $(B_3)$ and $(B_4)$ that two $5$-vertices are not adjacent to any $4$-vertex, then $\tilde{c}(\bar{f})\geq -2+2\times1=0$ by R4. For $(C_8)$, if three are three $5$-vertices lying on $\bar{f}$, it follows from Lemma \ref{configu1} $(B_6)$ that $\tilde{c}(\bar{f})\geq-2+2\times \frac{5}{6}+\frac{1}{2}=\frac{1}{6}>0$ by R3 and R4. Otherwise, $\tilde{c}(\bar{f})\geq-2+2\times\frac{2}{3}+\frac{4}{3}=\frac{2}{3}>0$ by R4-R6. For $(C_{9})$, $\tilde{c}(\bar{f})\geq-2+2\times \frac{2}{3}+2\times \frac{1}{2}=\frac{1}{3}>0$ by R3 and R4. For $(C_{10})$, $\tilde{c}(\bar{f})\geq-2+\frac{2}{3}+\frac{3}{2}=\frac{1}{6}>0$ by R4 and R5.

Let $v$ be a $k$-vertex of $G_1$. Suppose $k=4$, the initial charge remain unchanged. Suppose $k=5$. If $f_3(v)=0$, then there are at most five bad $4$-vertices, thus $\tilde{c}(v)\geq 1-5\times \frac{1}{6}=\frac{1}{6}>0$ by R1. If $f_3(v)=1$, we consider whether the triangle is bad. When the triangle is bad, we obtain that $\tilde{c}(v)\geq 1-\max\{\frac{2}{3}+\frac{1}{6},1\}=0$ by R1 and R3.1. Otherwise, $\tilde{c}(v)\geq 1-\max\{\frac{2}{3}+2\times \frac{1}{6}, \frac{1}{2}+3\times \frac{1}{6},1\}=0$ by R3.1 and R3.2. If $f_3(v)=2$, then the two triangles must form a doubleton $\bar{f}$ since $G_1$ is hopper-free, then we just need to consider the configurations in Figure \ref{doub} exclude $(C_1)$, $(C_2)$, $(C_3)$, $(C_6)$, then we obtain $\tilde{c}(v)\geq 1-\max\{1, \frac{5}{6}+\frac{1}{6}, \frac{2}{3}+2\times\frac{1}{6}\}=0$ by R4. Suppose $k=6$. If $f_3(v)=0$, then there are at most six bad $4$-vertices, thus $\tilde{c}(v)\geq 2-6\times \frac{1}{6}=1>0$ by R1. If $f_3(v)=1$, there are at most four bad $4$-vertices, thus $\tilde{c}(v)\geq 2-\frac{4}{3}-4\times \frac{1}{6}=0$ by R1 and R5.1. If $v$ is incident with two consecutive triangles $\bar{f}$, then we should consider $(C_6)$ and $(C_{10})$. As for $(C_6)$, if there are three $4$-vertices lying on
$\bar{f}$, then there exists no $4$-vertex among the remaining neighbors of $6$-vertex which are not on $\bar{f}$, thus $\tilde{c}(v)\geq 2-2=0$ by R5. Otherwise, $\tilde{c}(v)\geq 2-\frac{3}{2}-3\times \frac{1}{6}=0$ by R1 and R5. For $(C_{10})$, there are at most three bad $4$-vertices among the remaining neighbors of $6$-vertex which are not on $\bar{f}$ by Lemma \ref{configu1}, thus $\tilde{c}(v)\geq 2-\frac{3}{2}-3\times \frac{1}{6}=0$ by R5. Suppose $k\geq7$. If $f_3(v)=0$, then $\tilde{c}(v)\geq k-4-7\times \frac{1}{6}=\frac{6(k-7)+11}{6}$>0 by R1. If $f_3(v)=1$, then $\tilde{c}(v)\geq k-4-\frac{4}{3}-\frac{1}{6}(k-2)=\frac{5(k-7)+5}{6}>0$ by R1 and R6. If $f_3(v)=2$, then $\tilde{c}(v)\geq k-4-2-\frac{1}{6}(k-3)=\frac{5(k-7)+2}{6}>0$ by R6.

Hence, we complete the proof of Theorem \ref{thm1}.

\section{Proof of Theorem 2}\label{hfgra}
The notation we need in this section follows from Section \ref{nota}. Let $\mathcal{F}=\mathcal{H}_2$, and let $G_2$ be a counterexample to Theorem \ref{thm2} with minimum number of vertices. Fix a plane embedding of $G_2$, by minimality of $G_2$, we get that $G_2$ is connected. Let $L$ be a list assignment on $V(G_2)$ where each vertex receives at least five colors. Similarly, we have the following lemma to forbid some configurations in $G_2$.

\subsection{Reducible subgraphs}
\begin{figure}[H]
\begin{center}
\includegraphics[scale=0.8]{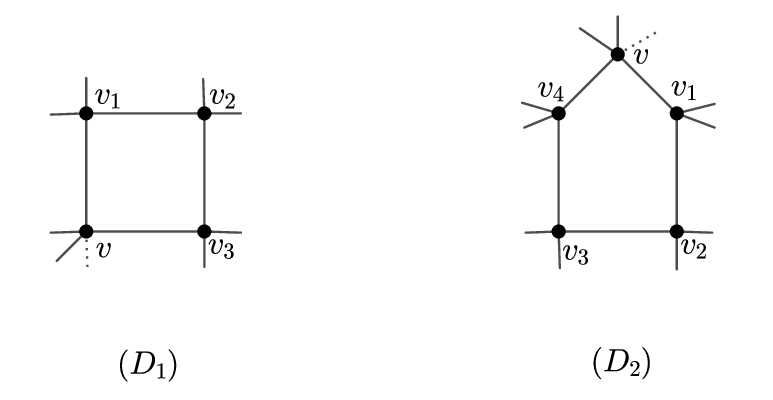}\\
\caption{Reducible graphs}
\label{hou}
\end{center}
\end{figure}

\begin{lemma}\label{CP}
Let $G_2\in \mathcal{H}_2$. If $G_2$ contains one of the following configurations (see Figure \ref{hou}),
\begin{description}
\item[$(D_1)$] A cycle $vv_1v_2v_3v$ such that $d(v)\leq5$ and $d(v_1)=d(v_2)=d(v_3)=4$;
\item[$(D_2)$] A cycle $vv_1v_2v_3v_4v$ such that $d(v)\leq5$ and $d(v_1)=d(v_2)=d(v_3)=d(v_4)=4$.
\end{description}
Then $G_2$ contains a $(\mathcal{H}_2,5)$-boundary-reducible induced subgraph with empty boundary.
\end{lemma}

\begin{proof}

In the following, we mainly consider both cases with $d(v)=5$.

To proof $(D_1)$. Let $H$ be the subgraph of $G_2$ induced by $\{ v, v_1, v_2, v_3\}$ and set the boundary $B=\emptyset$.

(\textbf{FIX}): Note that $v_i$ has three available colors for each $i\in \{1,2,3\}$, and $v$ has two available colors. If $v$ is fixed, then we greedily color $v$, $v_1$, $v_2$, $v_3$ in order. If $v_i$ is fixed for each $i\in\{1,2,3\}$, then we first color $v_i$ and then $H-v_i$ can be colored by Lemma \ref{degree}.

(\textbf{FORB}): Note that $|S|\leq 2$ by the definition of $S$.  When $|S|=2$, if $S=\{v,v_2\}$, then we give $v$ a preferred color and then greedily color $v_1$, $v_2$, $v_3$ in order. If $S=\{v_1, v_3\}$, then $H$ can be colored by Lemma \ref{degree}. When $|S|=1$, \textbf{(FORB)} will be implied by \textbf{(FIX)}.

To proof $(D_2)$, let $H$ be the subgraph of $G_2$ induced by $\{ v, v_1, v_2, v_3, v_4\}$, and set the boundary $B=\emptyset$.

(\textbf{FIX}): Note that $v_i$ has three available colors for each $i\in \{1,2,3,4\}$, and $v$ has two available colors. If $v$ is fixed, then we greedily color $v$, $v_1$, $v_2$, $v_3$, $v_4$ in order. If $v_i$ is fixed for each $i\in\{1,2,3,4\}$, then we first color $v_i$ and then $H-v_i$ can be colored by Lemma \ref{degree}.

(\textbf{FORB}): When $|S|=3$, then $S$ can be chosen as the following triples: $\{v,v_1,v_2\}$, $\{v,v_1,v_4\}$, $\{v,v_2,v_3\}$, $\{v,v_3,v_4\}$, $\{v_1,v_2,v_3\}$, $\{v_2,v_3,v_4\}$. Up to symmetry, we only discuss $S=\{v,v_1,v_2\}$. In this situation, we can greedily color $v$, $v_1$, $v_2$, $v_3$, $v_4$ in order. When $|S|=2$, assume that $S=\{v,v_1\}$, then we give $v$ a preferred color and then greedily color $v_1$, $v_2$, $v_3$, $v_4$ in order. The rest pairs can be handled in a similar fashion, (\textbf{FORB}) holds.
\end{proof}

\subsection{Discharging}
Note that $G_2$ does not contain any configurations shown in Lemma \ref{CP}. By \cite{MT}, $G_2$ also has no $3^-$-vertex. Since $G_2$ is also a plane graph, by Euler's Formula, we obtain

$$\sum_{v\in V(G_2)}(d(v)-6)+\sum_{f\in F(G_2)}(2d(f)-6)=-12. \eqno{(\mathrm{III})}$$

Now we define an initial weight function on $V(G_2)\cup F(G_2)$ by letting
\begin{align*}
\begin{split}
c(x)=\left\{
 \begin{array}{ll}
  d_{G_2}(x)-6    & $if$ \ x=v\in V(G_2), \\
  2d_{G_2}(x)-6    & $if$ \ x=f\in F(G_2).
 \end{array}
 \right.
 \end{split}
\end{align*}

Since the total sum of charges are the negative number $-12$, we shall now redistribute the charge, without changing its sum, such that the sum is nonnegative. This contradiction will prove the Theorem \ref{thm2}. Finally, we apply the following rules to redistribute the initial charge that leads to a new charge $\hat{c}$.

For brevity, a face $f$ is called \emph{bad} in $G_2$ if it is incident with exactly $(d(f)-1)$ $4$-vertices, otherwise it is \emph{good}.

\begin{description}
\item[$\bf R1.$] Each bad $4$-face sends $\frac{2}{3}$ to its incident $4$-vertex;
\item[$\bf R2.$] Each good $4$-face sends $\frac{1}{2}$ to its incident vertex $5^{-}$-vertex;
\item[$\bf R3.$] Each $5$-face sends $a$ to its incident vertex $v$,\\
$a=\begin{cases}
1 \ \ & $if$ \ d(v)=4, \\
\frac{1}{2} \ \ & $if$ \ d(v)=5.
\end{cases}$
\item[$\bf R4.$] Each $6^{+}$-face sends $1$ to its incident vertex $5^{-}$-vertex.
\end{description}

\begin{fact}
For each $v$ with $d(v)\geq4$, we have $f_3(v)\leq \lfloor\frac{2d(v)}{3}\rfloor$.
\end{fact}

Now we shall show that $\hat{c}(x)\geq0$ for all $x\in V(G_2)\cup F(G_2)$. Let $f$ be a face of $G_2$. If $d(f)=3$, we keep the initial charge. Suppose $d(f)=4$. If $f$ is bad, then there must be a $6^{+}$-vertex in $f$ by Lemma \ref{CP} $(D_1)$, thus $\hat{c}(f)\geq 2-3\times \frac{2}{3}=0$ by R1. Otherwise, $n_4(f)\leq2$, then $\hat{c}(f)\geq 2-4\times \frac{1}{2}=0$ by R2. Suppose $d(f)=5$. If $f$ is bad, then there must be a $6^{+}$-vertex in $f$ by Lemma \ref{CP} $(D_2)$, thus $\hat{c}(f)\geq 4-4\times1=0$ by R3. Otherwise, $n_4(f)\leq3$, then $\hat{c}(f)\geq 4-3\times1-2\times \frac{1}{2}=0$. If $d(f)\geq6$, then the number of $5^{-}$-vertices incident with $f$ is at most $d(f)$, thus $\hat{c}(f)\geq 2d(f)-6-d(f)=d(f)-6\geq0$ by R4.

Let $v$ be a $k$-vertex of $G_2$. Suppose $k=4$, then $f_3(v)\leq2$. If $f_3(v)=0$, then $\hat{c}(v)\geq -2+4\times \frac{1}{2}=0$ by R2. If $1\leq f_3(v)\leq2$, then $f_{5^{+}}(v)\geq2$, thus $\hat{c}(v)\geq -2+2\times1=0$ by R3. Suppose $k=5$, then $f_3(v)\leq3$. If $f_3(v)=0$, then $\hat{c}(v)\geq -1+5\times \frac{1}{2}>0$ by R2. If $1\leq f_3(v)\leq3$, then $f_{5^{+}}(v)\geq2$, thus $\hat{c}(v)\geq -1+2\times \frac{1}{2}=0$ by R3.
Suppose $k\geq6$, then $\hat{c}(v)=d(v)-6\geq0$.

Hence, for all $x\in V(G_2) \cup F(G_2)$, we have $\hat{c}(x)\geq0$. We complete the proof of Theorem \ref{thm2}.

\section{Acknowledgement}
The author would like to thank the the anonymous referees for their valuable remarks, which greatly improve the paper.

\end{document}